\documentclass[10pt,a4paper]{siamltex1213}
\usepackage[english]{babel}
\usepackage{amsmath}
\usepackage{amsfonts}
\usepackage{amssymb}
\usepackage{graphicx}
\usepackage{subcaption}
\usepackage{xcolor}
\usepackage{booktabs}
\usepackage{mdframed}

\usepackage{tikz}
\usetikzlibrary{shapes,arrows}

\author{Gijs L. Kooij\footnotemark[2] \and Mike A. Botchev\footnotemark[3] \and Bernard J. Geurts\footnotemark[2]\ \footnotemark[4]}
\title{A block Krylov subspace implementation of the time-parallel Paraexp method and its extension for nonlinear partial differential equations}

\newcommand{\mbf}[1]{\mathbf{#1}}

\newcommand{\Rr}{\mathbb{R}}

\newtheorem{remark}{Remark}

\begin{document}
\maketitle

\renewcommand{\thefootnote}{\fnsymbol{footnote}}

\footnotetext[2]{Multiscale Modeling and Simulation, Faculty of EEMCS, University of Twente, 7500 AE Enschede, The Netherlands.
Email: \email{g.l.kooij@utwente.nl} (first and corresponding author), \email{b.j.geurts@utwente.nl} (third author).}
\footnotetext[3]{Mathematics of Computational Science, Faculty of EEMCS, University of Twente, 7500 AE Enschede, The Netherlands (\email{m.a.botchev@utwente.nl}).}
\footnotetext[4]{Faculty of Applied Physics, Fluid Dynamics Laboratory, Eindhoven University of Technology, Eindhoven, The Netherlands.}

\renewcommand{\thefootnote}{\arabic{footnote}}

\slugger{mms}{xxxx}{xx}{x}{x--x}

\begin{abstract}
A parallel time integration method for nonlinear partial differential equations is proposed.  It is based on a new implementation of the Paraexp method for linear partial differential equations (PDEs) employing a block Krylov subspace method. For nonlinear PDEs the algorithm is based on our Paraexp implementation within a waveform relaxation. The initial value problem is solved iteratively on a complete time interval. Nonlinear terms are treated as a source term, provided by the solution from the previous iteration. At each iteration, the problem is decoupled into independent subproblems by the principle of superposition. The decoupled subproblems are solved fast by exponential integration, based on a block Krylov method. The new time integration is demonstrated for the one-dimensional advection-diffusion equation and the viscous Burgers equation. Numerical experiments confirm excellent parallel scaling for the linear advection-diffusion problem, and good scaling in case the nonlinear Burgers equation is simulated.
\end{abstract}



\pagestyle{myheadings}
\thispagestyle{plain}
\markboth{KOOIJ ET AL.}{A TIME-PARALLEL EXPONENTIAL INTEGRATOR}

\section{Introduction}

Recent developments in hardware architecture, bringing to practice computers with hundreds of thousands of cores, urge the creation of new, as well as a major revision of existing numerical algorithms~\cite{Dongarra:2011aa}. To be efficient on such massively parallel platforms, the algorithms need to employ all possible means to parallelize the computations. When solving partial differential equations (PDEs) with a time-dependent solution, an important way to parallelize the computations is, next to the parallelization across space, parallelization across time. This adds a new dimension of parallelism with which the simulations can be implemented. In this paper we present a new time-parallel integration method extending the Paraexp method~\cite{Gander:2013aa} to nonlinear partial differential equations using Krylov methods and waveform relaxation.

Several approaches to parallelize the simulation of time-dependent 
solutions in time can be distinguished. The first important class of the methods are the waveform relaxation methods~\cite{Botchev:2014aa,Miekkala:1987aa,Newton:1983aa,White:1985aa}, including the space-time multigrid methods for 
parabolic PDEs~\cite{Burrage:1995aa,Hackbusch:1984aa,Horton:1995aa,Lubich:1987aa}. The key idea is to start with an approximation to the 
numerical solution for the whole time interval of interest and 
update the solution, solving an easier-to-solve approximate 
system in time. The Parareal method~\cite{Lions:2001aa}, which attracted  significant attention recently, is a prime example related to the class of waveform relaxation methods~\cite{Gander:2007aa}.

Parallel Runge--Kutta methods and general linear methods, where the parallelism is determined and restricted by the number of stages or steps, form another class of the time-parallel methods~\cite{Burrage:1995aa,de-Swart:1997aa,Houwen:1998aa}. Time-parallel schemes can also be obtained by treating the time as an additional space dimension and solving for values at all time steps at once~\cite{Dolgov:2012aa,Zhu:2000aa}. This approach requires significantly more memory and is used, e.g., in multidimensional tensor computations and (discontinuous) Galerkin finite element methods~\cite{Cella:1975aa,Hughes:1988aa,Lions:1961aa,Vegt:2002aa}. Recently a ``parallel full approximation scheme in space and time'' (PFASST) was introduced, which is based on multigrid methods~\cite{Emmett:2012aa}. PFASST was observed to have a generally improved parallel efficiency compared to Parareal. Also, we mention parallel methods that facilitate parallelism by replacing the exact solves in implicit schemes by approximate ones~\cite{Botchev:2001aa}, and parallel methods based on operator splitting~\cite{Botchev:2004aa}. Finally, there is the Paraexp method~\cite{Gander:2013aa} for parallel integration of linear initial-value problems. This algorithm is based on the fact that linear initial-value problems can be decomposed into homogeneous and nonhomogeneous subproblems. The homogeneous subproblems can be solved fast by an exponential integrator, while the nonhomogeneous subproblems can be solved a traditional time-stepping method.

In this paper we propose a time parallel method which is based on the matrix exponential time stepping technique and waveform relaxation. Our method adopts the Paraexp method within a waveform relaxation framework, which enables the extension of parallel time integration to nonlinear PDEs. The method is inspired by and relies on recent techniques in Krylov subspace matrix exponential methods such as shift-and-invert acceleration~\cite{Gockler:2014aa,Moret:2004aa,Eshof:2006aa}, restarting~\cite{Afanasjew:2008aa,Celledoni:1997aa,Eiermann:2011aa,Niehoff:2006aa,Tal-Ezer:2007aa} and using block Krylov subspaces~\cite{Botchev:2013aa} to handle nonautonomous PDEs efficiently. {\color{black}
The method also relies on a singular value decomposition (SVD) of  source terms in the PDE, which is used to construct the block Krylov subspace. To improve the efficiency of the method, the SVD is truncated by retaining only the relatively large singular values. We show in theory and in practice that for a source term, that can be approximated by a smooth function, the singular values decay algebraically with time interval size. In that case, a truncated SVD approximation of the source term is adequate.
}

The contribution of this paper is more specifically as follows. First, to solve systems of linear ODEs, we propose an implementation of the Paraexp method, based on the exponential block Krylov (EBK) method, introduced in~\cite{Botchev:2013aa}. {The EBK method is a Krylov method that approximates the exact solution of a system of nonhomogeneous linear ordinary differential equations by a projection onto a block Krylov subspace.} Our Paraexp implementation does not involve a splitting of homogeneous and nonhomogeneous subproblems, which leads to a better parallel effiency. Second, we extend our EBK-based implementation of the Paraexp method to nonlinear initial-value problems. We solve the problem iteratively on a chosen time interval. In our case, the nonlinear term of the PDE is incorporated as a source term, which is updated after every iteration with the latest solution. Third, we show that our Paraexp-EBK (PEBK) implementation has promising properties for time-parallelization. The PEBK method can be seen as a special continuous-time waveform relaxation method, which is typically more efficient than standard waveform relaxation methods~\cite{Botchev:2014aa}. The PEBK method is tested for the one-dimensional advection-diffusion equation and the viscous Burgers equation in an extensive series of numerical experiments. Since we are interested in time-parallel solvers for large-scale applications in computational fluid dynamics (CFD), such as turbulent flow simulations, we also test our method with respect to the diffusion/advection ratio and grid resolution, reflecting high Reynolds-number conditions. For example, the parallel efficiency of the Parareal algorithm was found to decrease considerably with increasing number of processors for the advection equation~\cite{Gander:2008aa,Ruprecht:2012aa,Staff:2005aa}. In contrast, the EBK method was found to have excellent weak scaling for linear problems. 

	The paper is organized as follows. In Section~\ref{sec:expIntegration} we describe the exponential time integration scheme and its parallelization. In Section~\ref{sec:ade}, we present numerical experiments with the advection-diffusion equation, and with the viscous Burgers equation in Section~\ref{sec:burgers}. Finally, a discussion and conclusions are outlined in Section~\ref{sec:conclusions}.

\section{Exponential time-integration}
\label{sec:expIntegration}

Our time integration method for linear and nonlinear PDEs is based on the exponential block Krylov (EBK) method~\cite{Botchev:2013aa}. In this section, we first provide a brief description of the EBK method. Then, we extend the EBK method to integrate nonlinear PDEs in an iterative way. Finally, the parallelization of the time integrator is discussed.

\subsection{Exponential block Krylov method}
The EBK method is a time integrator for linear systems of nonhomogeneous ordinary differential equations (ODEs). Details of the method are given in~\cite{Botchev:2013aa}. We follow the \emph{method of lines} approach~\cite{Verwer:1984aa}, i.e., the PDE is discretized in space first. We start with linear PDEs. After applying a spatial discretization to the PDE, we obtain the initial value problem,
\begin{equation}
\begin{aligned}
	\mbf{u}'(t) &= A\mbf{u}(t) + \mbf{g}(t), \\
	\mbf{u}(0) &= \mbf{u}_0,
	\label{eq:ivp2}
\end{aligned}
\end{equation}
where $\mbf{u}(t)$ is a vector function $\mbf{u}(t): \mathbb{R} \rightarrow \mathbb{R}^n$, $A \in \mathbb{R}^{n \times n}$ is a square matrix, and $\mbf{g}(t) \in \mathbb{R}^n$ a time-dependent source term. The vector ${\mbf{u}}$ can readily be identified with the solution to the PDE in terms of its values taken on a computational grid in physical space. In Section~\ref{sec:ebk_nonlinear}, we will introduce a more general term $\mbf{g}(t,\mbf{u}(t))$, which contains the nonlinear terms of the PDE. The matrix $A$ is typically large and sparse, depending on the spatial discretization method. The dimension of the system, $n$, corresponds to the number of degrees of freedom used in the spatial discretization. Exponential integrators approximate the exact solution of the semi-discrete system~\eqref{eq:ivp2}, formulated in terms of the matrix exponential. 

The first step in the EBK method is to approximate $\mbf{g}(t)$ with a function that can be treated easier. Here, $\mbf{g}(t)$ is approximated based on a truncated singular value decomposition (SVD)~\cite{Botchev:2009aa} as follows. The source term is sampled at $s$ points in time, $0 = t_1 < t_2 < \ldots < t_{s-1} < t_s = \Delta T$, in the integration interval $[0,\Delta T]$, and the  source samples form the matrix
\begin{equation}
G = \left[\,\mbf{g}(t_1)\quad\mbf{g}(t_2)\quad\ldots\quad\mbf{g}(t_s)\,\right] \in \mathbb{R}^{n \times s}.
\end{equation}
We make a natural assumption that the typical number of time samples $s$, necessary for a good approximation of $\mbf{g}(t)$, is much lower than $n$: $s \ll n$. 
Since we assume $s \ll n$, it is more economical to calculate the so-called thin SVD of $G$~\cite{Golub:1996aa}, instead of the full SVD, without loss of accuracy. In this case, the thin SVD of $G$ is
\begin{equation}
G = \tilde{U} \tilde{\Sigma} \tilde{V}^T,
\end{equation}
where $\tilde{\Sigma} \in \mathbb{R}^{s \times s}$ is a diagonal matrix containing the singular values $\sigma_1 \geq \sigma_2 \geq \ldots \geq \sigma_s$, and $\tilde{U} \in \mathbb{R}^{n \times s}$, $\tilde{V} \in \mathbb{R}^{s \times s}$ are matrices with orthonormal columns. The thin SVD can be approximated by a truncation in which the $m < s$ largest singular values are retained. As seen from this truncated SVD, the samples of $\mbf{g}(t)$ can be approximated by linear combinations of the first $m$~columns of $\tilde{U}$, i.e.,
\begin{equation}
	\mbf{g}(t) \approx U \mbf{p}(t).
	\label{eq:svd_appr}
\end{equation}
where $U \in \mathbb{R}^{n \times m}$ is formed by the first $m$ columns of $\tilde{U}$, $\mbf{p}(t) \in \mathbb{R}^m$ is obtained by an interpolation of the coefficients in these linear combinations~\cite{Botchev:2009aa}. There are several possible choices for $\mbf{p}(t)$, among which, cubic piecewise polynomials. Then, the approximation error in the source term, $\Vert \mbf{g}(t) - U\mbf{p}(t) \Vert$, can be easily controlled within a desired tolerance, depending on the number of samples in $[0,\Delta T]$, and the number of singular values truncated (see \cite{Botchev:2013aa}).

{\color{black} The number of retained singular values, $m$, is equal to the block width in the block Krylov subspace. The efficiency of the EBK method therefore depends on how many singular values are required for an accurate approximation of the source term~\eqref{eq:svd_appr}. In applications of the method this parameter $m$ can be varied and practical convergence can be assessed.} If the subinterval $\Delta T$ is small and the source function $\mbf{g}(t)$ can be well approximated by a smooth function, then the singular values of the samples of $\mbf{g}$ decrease rapidly, thus allowing for an efficient low rank approximation in~\eqref{eq:svd_appr}. Furthermore, let the subinterval length $\Delta T$
be small in the sense that $(\Delta T)^2\ll\Delta T$ (in dimensionless time units). Of course, the assumption that $\Delta T$ is small, is not always realistic and we comment on how it can be partially relaxed below in Remark~\ref{Rem1}.

Denote by  $\mbf{g}^{(j)}$ the $j$th derivative of $\mbf{g}$ and assume that
the constants
$$
M_j :=\max_{\xi\in[t,t+\Delta T]} \| g^{(j)}(\xi)\|
\qquad
j = 0,\dots , n,
$$
are bounded. 
We can then show that the singular values of the samples decrease, starting from a thin QR~factorization of $G$.

Consider the thin QR~factorization of $G$, $G=QR$, where $Q\in\Rr^{n\times s}$
  has orthonormal columns and $R\in\Rr^{s\times s}$ is an upper triangular
  matrix with nonnegative diagonal entries. As Theorem~1 in~\cite{Hochbruck:2008aa} states, the entries $r_{ik}$ of $R$ satisfy
  $$
  |r_{jk}|\leq C_j M_{j-1}(\Delta T)^{j-1}, \quad k\geq j,
  $$
  and $r_{jk}=0$ for $k<j$ because $R$ is upper triangular.
  In this estimate the constants $C_j$ depend only on the
  points $t_1$, \dots, $t_s$ at which the samples are computed and
  do not depend on $\mbf{g}$ and $\Delta T$.

  Since $G$ and $R$ have the same singular values, we consider now the
  singular values of $R$.  According to~\cite[3.1.3]{Horn:1991aa},
  \begin{equation}
  \label{estim1}  
  \sigma_{j+1}\geq \|R_j\|_2,\qquad j=1,\dots,s-1,
  \end{equation}
  where $R_j$ is a matrix obtained from $R$ by skipping
  $j$ rows or columns, chosen arbitrarily.  Since the entries
  in $R$ decrease rowwise as $|r_{jk}|=O(\Delta T)^{j-1}$, to have
  the sharpest estimate in~\eqref{estim1}, we choose $R_j$ to be
  the matrix $R$ with the first $j$ rows skipped.
  To bound the 2-norm of $R_j$ we use~\cite[Chapter~5.6]{Horn:1986aa}
  $$
  \|R_j\|_2 \leq \sqrt{\|R_j\|_1\|R_j\|_{\infty}}
  $$
  and note that for $j=1,\dots,s-1$
  $$
  \|R_j\|_{\infty} \leq (s-j)C_jM_{j-1}(\Delta T)^{j-1},
  \qquad
  \|R_j\|_1 = O(\Delta T)^{j-1}.
  $$
  Thus, we obtain $\|R_j\|_2 = \sqrt{s-j}\,O(\Delta T)^{j-1}$, which,
  together with~\eqref{estim1}, yields the following result.

\begin{theorem}
  \label{The1}
  Let $\mbf{g}(t): \Rr\rightarrow\Rr^n$ be a smooth function
  such that the constants
  $$
  M_j :=\max_{t\in[0,\Delta T]} \| g^{(j)}(t)\|
  \qquad
  j = 0,\dots , n,
  $$
  are bounded.  Furthermore, let the subinterval length $\Delta T$
  be small in the sense that $(\Delta T)^2\ll\Delta T$.  Then for the
  singular values $\sigma_j$, $j=1,\dots,s$, of the sample matrix
  $$
  G = 
  \begin{bmatrix}
  \mbf{g}(t_1) & \mbf{g}(t_2) & \dots & \mbf{g}(t_s)  
  \end{bmatrix}\in \Rr^{n\times s}
  $$
  holds
  \begin{equation}
    \label{sigma_bound}
    \sigma_{j+1}= \sqrt{s-j}\, O(\Delta T)^j, \qquad j=1,\dots,s-1.
  \end{equation}
\end{theorem}

\begin{remark}
  \label{Rem1}
  We may extend the result of Theorem~\ref{The1} to the case
  of a large $\Delta T$ if the constants $M_j$ are bounded
  more strongly or even decay with $j$.
  Indeed, if $\Delta T$ is large, we can consider the function
  $\tilde{\mbf{g}}(\tilde{\xi}):=\mbf{g}(t+\xi q)$, with $q$ chosen such
  that $\tau:=\Delta T/ q$ is small.  As the function $\tilde{\mbf{g}}$
  takes the same values for $\tilde{\xi}\in [0,\tau]$ as
  the function $\mbf{g}$ for $\xi\in[t,t+\Delta T]$,
  Theorem~\ref{The1}
  formally holds for $\tilde{\mbf{g}}$ with $\Delta T$ replaced by $\tau$
  and $M_j$ multiplied by $q^j$.  The coefficients
  in the $O$~symbol in~\eqref{sigma_bound} may now grow with
  the powers of $\Delta T$, thus rendering the result meaningless
  unless a stricter assumption on $M_j$ is made.
\end{remark}

The truncated SVD approximation of the source term~\eqref{eq:svd_appr} facilitates the solution of the initial value problem (IVP) in~\eqref{eq:ivp2} by the block Krylov subspace method~\cite{Botchev:2013aa}. Here, the block Krylov subspace at iteration $l$ is defined as
\begin{equation}
	\mathcal{K}_l (A,U) := \text{span} \left\{ U, AU, A^2U, \ldots, A^{l-1}U \right\}.
	\label{eq:blockKrylovSubspace}
\end{equation}
Furthermore, the block Krylov subspace method can be accelerated by the shift-and-invert technique~\cite{Gockler:2014aa,Moret:2004aa,Eshof:2006aa} and, to keep the Krylov subspace dimension bounded, implemented with restarting~\cite{Botchev:2013ab,Celledoni:1997aa}.

\subsection{Parallelization of linear problems}
\label{sec:parlinear}
\begin{figure}
\includegraphics[scale=1]{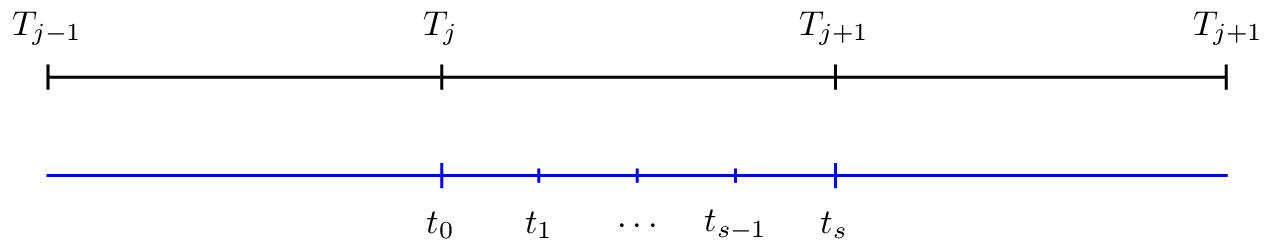}
\caption{Time samples of the source term on the subinterval $[T_j,T_{j+1}]$.}
\label{fig:subgrid}
\end{figure}

In this section, the parallelization of the linear ODE solution is discussed. As we will see, the method is equivalent to the Paraexp method, but there are differences in implementation leading to a better parallel efficiency (see Section~\ref{sec:ade_efficiency}). Linear IVPs can be solved parallel in time by using the principle of superposition. First, the time intervals is divided into $P$ non-overlapping subintervals,
\[
	0 = T_0 < T_1 < \ldots < T_P = T,
\]
where $P$ is also the number of processors that can be used in parallel. We introduce the shift of $\mbf{u}(t)$ with respect to the initial condition, $\mbf{u}(t) = \mbf{u}_0 + \widehat{\mbf{u}}(t)$. The shifted variable $\widehat{\mbf{u}}(t)$ solves the IVP with homogeneous initial conditions,
\begin{equation}
\begin{aligned}
	\widehat{\mbf{u}}'(t) &= A \widehat{\mbf{u}}(t) + \widehat{\mbf{g}}(t), \quad t \in [0,T],\\
	\widehat{\mbf{u}}(0) &= \mbf{0},
	\label{eq:linivp_shift}
\end{aligned}
\end{equation}
where
\begin{equation}
\widehat{\mbf{g}}(t) := A\mbf{u}_0 + \mbf{g}(t).
\end{equation}
The source term is approximated using $s$ time samples per subinterval, as illustrated in Fig.~\ref{fig:subgrid}.
The shifted IVP~\eqref{eq:linivp_shift} can be decoupled into independent subproblems by the principle of superposition. To each subinterval $[T_j,T_{j-1}]$ we associate a part of the source term $\widehat{\mbf{g}}(t)$, defined for $j = 1,\ldots,P$ as
\begin{equation}
	\widehat{\mbf{g}}_{j}(t) = \left\{
		\begin{array}{ll}
			\widehat{\mbf{g}}(t), & \text{for } T_{j-1} \leq t < T_{j},\\
			0, & \text{otherwise} .
		\end{array}
	\right.
	\label{eq:piecewiseSource_linivp}
\end{equation}
The expected parallel speedup is based on the observation that the solution to \eqref{eq:linivp_shift} is given by a variation-of-constant formula (see, e.g., \cite{Hochbruck:2010aa}),
\begin{equation}
\widehat{\mbf{u}}(T) = \int_0^T \exp \left( (T - s)A \right) \widehat{\mbf{g}}(s)\,ds = \sum_{j=1}^P \int_{T_{j-1}}^{T} \exp \left( (T - s)A \right) \widehat{\mbf{g}}_{j}(s)\,ds,
\end{equation}
where the integrals $\int_{T_{j-1}}^{T_j} \ldots ds$ can be evaluated independently and in parallel. More precisely, by exploiting the linearity of the problem we decompose the IVP~\eqref{eq:ivp_shift} into $P$ independent subproblems,
\begin{equation}
\begin{aligned}
	\mbf{v}_{j}'(t) &= A \mbf{v}_{j}(t) + \widehat{\mbf{g}}_{j}(t), \quad t \in [0,T],\\
	\mbf{v}_{j}(0) &= \mbf{0}.
	\end{aligned}
	\label{eq:linivpSubproblem}
\end{equation}
where $\mbf{v}_j(t)$ is referred to as a subsolution of the total problem. The subproblems are independent and can be solved in parallel. We integrate the subproblem individually with the EBK method. Observe that the source is only nonzero on the subinterval $[T_{j-1},T_j)$. We follow the ideas of the Paraexp method~\cite{Gander:2013aa} and note that the nonhomogeneous part of the ODE requires most of the computational work in the EBK method. An accurate SVD approximation of the source term~\eqref{eq:svd_appr} generally requires more singular values to be retained, increasing the dimensions of the block Krylov subspace~\eqref{eq:blockKrylovSubspace}. According to the principle of superposition, the solution of the original problem is then given as
\begin{equation}
\mbf{u}(t) = \mbf{u}_0 + \sum_{j=1}^P \mbf{v}_{j}(t).
\label{eq:superposition}
\end{equation}
The summation of the subsolutions, $\mbf{v}_j(t)$ is the only communication required between the parallel processes. The parallel algorithm is summarized in Fig.~\ref{fig:algorithm_linear}. In principle, this algorithm is identical to the Paraexp method~\cite{Gander:2013aa}. The only practical difference is that in our implementation both the nonhomogeneous and the homogeneous part of the subproblems are solved by the EBK method. The original version of the Paraexp method assumes a convential time integration method to solve the ``difficult'' nonhomogeneous part, and a Krylov subspace method for exponential integration of the homogeneous part. Our implementation of the Paraexp method is compared numerically with the original one in Section~\ref{sec:ade_efficiency}.

Finally, additional parallelism could be exploited within the block Krylov subspace method itself. If the block size is $m$, then the $m$ matrix-vector products can be executed entirely in parallel, see~\cite{Li:1997aa,Nikishin:1995aa}. This approach could be applied in combination with the parallelization  described in this section.

\begin{figure}[hbtp]
\begin{mdframed}
\textsc{Algorithm parallel time integration.\\}
Given: $A$, $\mbf{u}_0$, $\mbf{g}(t)$, ...\\
Solve: $\mbf{u}'(t) = A\mbf{u}(t) + \mbf{g}(t), \quad \mbf{u}(0) = \mbf{u}_0.$\vspace*{5pt}\\
\hspace*{0em}\textbf{for} $j = 1,\ldots,P$ (in parallel)\\
\hspace*{1em}Solve nonhomogeneous part of the ODE~\eqref{eq:linivpSubproblem}, for $t\in [T_{j-1},T_j]$.\\
\hspace*{1em}Solve homogeneous part of the ODE~\eqref{eq:linivpSubproblem}, for $t \in [T_j,T_P].$\\
\hspace*{0em}\textbf{end for}\\
\hspace*{0em}Construct solution $\mbf{u}(t)$, see~\eqref{eq:superposition}.
\end{mdframed}
\caption{The algorithm for solving linear ODE systems with the Paraexp exponential block Krylov (PEBK) method.}
\label{fig:algorithm_linear}
\end{figure}

\subsection{Treatment of nonlinearities}
\label{sec:ebk_nonlinear}
The EBK method, designed to solve a linear system of nonhomogeneous ODEs~\eqref{eq:ivp2}, can be extended to handle nonlinear systems of ODEs by including the nonlinearities in the source term. The system is then solved iteratively, in such a way that the iterand $\mbf{u}_k(t)$ is updated on the entire interval $[0,T]$:
\begin{equation}
\mbf{u}_k(t) \mapsto \mbf{u}_{k+1}(t),
\end{equation}
where $k$ denotes the iteration index. We proceed in a few steps. First, the problem is reduced to the form of Eq.~\eqref{eq:ivp2}. The nonlinear IVP is therefore approximated by evaluating the source term with the current iterand, $\mbf{u}_k(t)$. Because the current iterand $\mbf{u}_k(t)$ is a known function, we can write the source term as an explicit function of time,
\begin{equation}
\mbf{g}_k(t) := \mbf{g}(t,\mbf{u}_k(t)).
\label{eq:source}
\end{equation}
The resulting initial value problem is,
\begin{equation}
\begin{aligned}
	\mbf{u}_{k+1}'(t) &= A\mbf{u}_{k+1}(t) + \mbf{g}_k(t), \\
	\mbf{u}_{k+1}(0) &= \mbf{u}_0.
	\label{eq:ivp3}
\end{aligned}
\end{equation}
This system is solved by the EBK method to achieve one iteration. This process is continued until the solution is sufficiently converged. This approach is similar to applying Picard or fixed-point iterations~\cite{Nevanlinna:1989aa} on the nonlinear term.
The source term can be further improved for nonlinear behaviour by using the Jacobian matrix, similar to a Newton--Rhapson method. In this case, we introduce the average Jacobian matrix,
\begin{equation}
J_k := \frac{1}{T} \int_0^T \left[ \frac{\partial{\mbf{g}_k}}{\partial u_1}\: \ldots \: \frac{\partial{\mbf{g}_k}}{\partial u_n} \right]\,dt ,
\label{eq:jacobian}
\end{equation}
which is averaged over the interval of integration $[0,T]$. Generally, the solution varies in time and so does the Jacobian matrix. The nonlinear remainder, with respect to the time-averaged state, is contained in the source term. Using this correction, we then find the recursive relation,
\begin{equation}
\begin{aligned}
	\mbf{u}_{k+1}'(t) &= [A + J_k] \mbf{u}_{k+1}(t) + \mbf{g}_k(t) - J_k \mbf{u}_{k}(t),\\
	\mbf{u}_{k+1}(0) &= \mbf{u}_0.
	\label{eq:ivp4}
\end{aligned}
\end{equation}
When converged, $\mbf{u}_{k+1}(t) = \mbf{u}_k(t)$, the terms containing $J_k$ eventually disappear.

\begin{remark}
\label{Rem2}
The iteration~\eqref{eq:ivp4} is well known in the literature on  
waveform relaxation methods~\cite{Lubich:1987aa,Miekkala:1987aa}.
Its convergence is given, e.g., in~\cite{Botchev:2013ab,Nevanlinna:1989aa} and, in case of an inexact iteration, in~\cite{Botchev:2014aa}. These results, in particular, show that the iteration~\eqref{eq:ivp4} converges superlinearly on \emph{finite} time intervals if the norm
of the matrix exponential of $t(A+J_k)$ decays exponentially
with $t$, i.e., the eigenvalues of $A+J_k$ lie in the left half-plane.
\end{remark}

\subsection{Parallelization of nonlinear problems}
\label{sec:parnonlinear}
Nonlinear IVPs are solved in an iterative way with the PEBK method. We follow the \emph{waveform relaxation} approach \cite{Botchev:2014aa,Lelarasmee:1982aa,Nevanlinna:1989aa}, that is, the problem is solved iteratively on the entire time interval of interest, $[0,T]$. The original problem is decomposed into a number of independent, ``easier'' subproblems on $[0,T]$ that can be solved iteratively, and in parallel. 

At each iteration the nonlinear term is treated as source term, which gives us a linear system of ODEs that can be integrated in parallel, see Section~\ref{sec:parlinear}. In the case of nonlinear problems, we take the average of the Jacobian matrix on each subinterval, which gives us the piecewise constant function for $j = 1,\ldots,P$:
\begin{equation}
J_k(t) := \frac{1}{T_{j} - T_{j-1}} \int_{T_{j-1}}^{T_j} \left[ \frac{\partial{\mbf{g}_k}}{\partial u_1}\: \ldots \: \frac{\partial{\mbf{g}_k}}{\partial u_n} \right]\,dt, \quad t \in [T_{j-1},T_j].
\label{eq:piecewiseJac}
\end{equation}
The IVP is then shifted with respect to the initial condition,
\begin{equation}
\begin{aligned}
	\widehat{\mbf{u}}_{k+1}'(t) &= [A + J_k(t)] \widehat{\mbf{u}}_{k+1}(t) + \widehat{\mbf{g}}_k(t), \quad t \in [0,T],\\
	\widehat{\mbf{u}}_{k+1}(0) &= \mbf{0},
	\label{eq:ivp_shift}
\end{aligned}
\end{equation}
where,
\begin{equation}
\widehat{\mbf{g}}_k(t) := A\mbf{u}_0 + \mbf{g}_k(t) + J_k(t) [\mbf{u}_0 - \widehat{\mbf{u}}_{k}(t)].
\label{eq:sourceParallel}
\end{equation}
The shifted IVP~\eqref{eq:ivp_shift} can be then solved in parallel. We follow an approach similar to the Paraexp method, as explained in Section~\ref{sec:parlinear}. Here, we define
\begin{equation}
	\widehat{\mbf{g}}_{j,k}(t) := \left\{
		\begin{array}{ll}
			\widehat{\mbf{g}}_k(t), & \text{for } T_{j-1} \leq t < T_{j},\\
			0, &  \text{otherwise} .
		\end{array}
	\right.
	\label{eq:piecewiseSource}
\end{equation}
Now, the problem is decomposed into $P$ independent subproblems.  The subsolutions ${\mbf{v}}_{j}$ result from convergence of ${\bf{v}}_{j,k}$ obeying:
\begin{align}
	\mbf{v}_{j,k+1}'(t) &= [A + J_k(t)] \mbf{v}_{j,k+1}(t) + \widehat{\mbf{g}}_{j,k}(t), \quad t \in [0,T],\\
	\mbf{v}_{j,k+1}(0) &= \mbf{0},
	\label{eq:ivpSubproblem}
\end{align}
which can be integrated with the EBK method. Here, we have introduced a second subindex $k$ to indicate the $j$th subsolution at the $k$th iteration. The total solution can then be updated according to the principle of superposition,
\begin{equation}
\mbf{u}_{k+1}(t) = \mbf{u}_0 + \sum_{j=1}^P \mbf{v}_{j,k+1}(t).
\label{eq:updateSolution}
\end{equation}
After each iteration, the solution is assembled, and the time-averaged Jacobian matrix is updated at each subinterval, according to Eq.~\eqref{eq:piecewiseJac}. Note that for linear systems of ODEs, no iterations are required at all. The parallel algorithm is summarized in Fig.~\ref{fig:algorithm}.

\begin{figure}[hbtp]
\centering
\tikzstyle{decision} = [ellipse, draw, fill=blue!20, 
    text width=4.5em, text badly centered, node distance=1.5cm, inner sep=0pt, height=2em]
\tikzstyle{block} = [rectangle, draw, fill=blue!20, 
    text width=8em, node distance=1.5cm, text centered, rounded corners, minimum height=2em]
\tikzstyle{line} = [draw, -latex', semithick]
\tikzstyle{cloud} = [draw, ellipse,fill=red!20, node distance=3cm,
    minimum height=2em]
\begin{tikzpicture}
    \node [block] (init) {Initialize $\mathbf{u}_0(t)$ on $[T_j,T_{j+1}]$};
    \node [block, below of=init, text width=10em] (update) {Update source term~\eqref{eq:source} and Jacobian matrix~\eqref{eq:jacobian}};
    \node [block, below of=update] (solve) {Solve IVP~\eqref{eq:ivp4} for $\mbf{u}_{k+1}(t)$};
    \node [block, below of=solve] (stop) {$k < K$?};
    \node [block, below of=stop] (next) {Next time interval \dots};
    \node [block, left of=solve, node distance=3.5cm] (increment) {$k + 1$};

	\path [line] (init) -- (update);
	\path [line] (update) -- (solve);
	\path [line] (solve) -- (stop);
	\path [line] (stop) -| node [auto,near start] {yes} (increment);
	\path [line] (increment) |- (update);
	\path [line] (stop) -- node[auto] {no} (next);
\end{tikzpicture}
\caption{Flow diagram of the serial EBK method for nonlinear PDEs. The algorithm stops after $K$ iterations, after which the solution is assumed to be converged.}
\label{fig:flowSerial}
\end{figure}
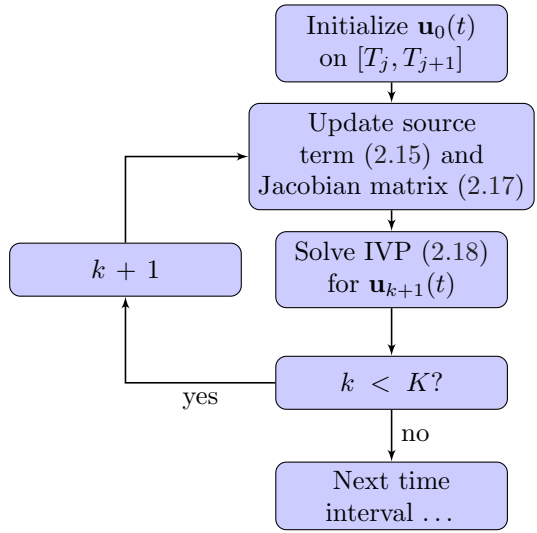

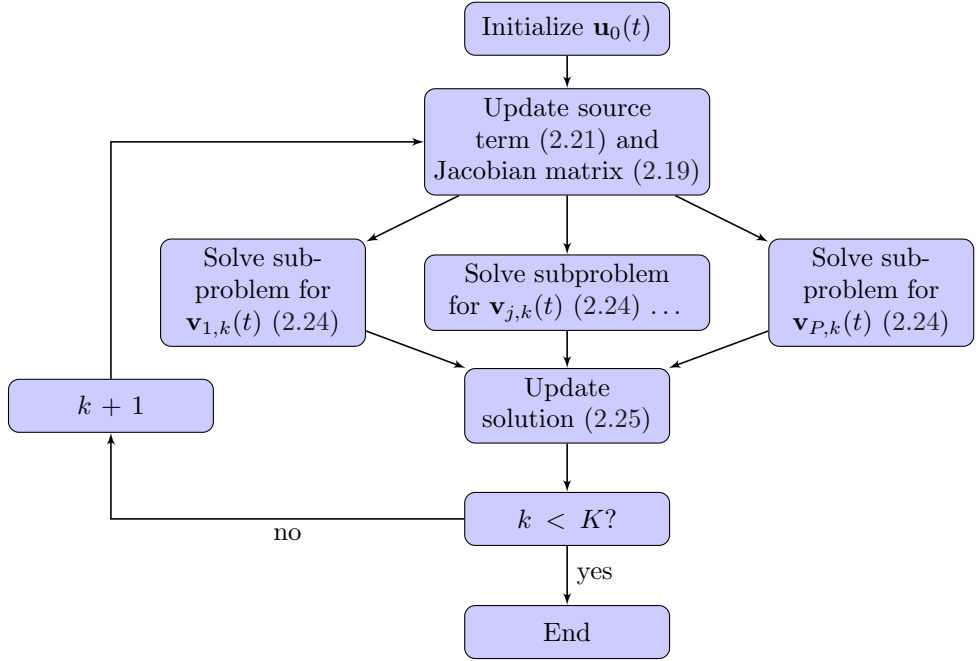
\begin{figure}[hbtp]
\centering
\tikzstyle{decision} = [ellipse, draw, fill=blue!20, 
    text width=4.5em, text badly centered, node distance=1.5cm, inner sep=0pt, height=2em]
\tikzstyle{block} = [rectangle, draw, fill=blue!20, 
    text width=7em, node distance=1.5cm, text centered, rounded corners, minimum height=2em]
\tikzstyle{line} = [draw, -latex', semithick]
\tikzstyle{cloud} = [draw, ellipse,fill=red!20, node distance=3cm,
    minimum height=2em]
\begin{tikzpicture}
    \node [block] (init) {Initialize $\mathbf{u}_0(t)$};
    \node [block, below of=init, text width=10em] (update) {Update source term~\eqref{eq:sourceParallel} and Jacobian matrix~\eqref{eq:piecewiseJac}};
    \node [block, below of=update, node distance=2cm, text width=10em] (sub2) {Solve subproblem for $\mathbf{v}_{j,k}(t)$~\eqref{eq:ivpSubproblem} \dots};
    \node [block, left of=sub2, node distance=4cm] (sub1) {Solve subproblem for $\mathbf{v}_{1,k}(t)$~\eqref{eq:ivpSubproblem}};
    \node [block, right of=sub2, node distance=4cm] (sub3) {Solve subproblem for $\mathbf{v}_{P,k}(t)$~\eqref{eq:ivpSubproblem}};
    \node [block, below of=sub2] (update2) {Update solution~\eqref{eq:updateSolution}};
    \node [block, below of=update2] (stop) {$k < K$?};
    \node [block, below of=stop] (end) {End};
    \node [block, left of=update2, node distance=6cm] (increment) {$k+1$};

\path [line] (init) -- (update);
\path [line] (update) -- (sub1);
\path [line] (update) -- (sub2);
\path [line] (update) -- (sub3);
\path [line] (sub1) -- (update2);
\path [line] (sub2) -- (update2);
\path [line] (sub3) -- (update2);
\path [line] (update2) -- (stop);
\path [line] (stop) -- node[auto] {yes} (end);
\path [line] (stop) -| node[auto,near start] {no} (increment);
\path [line] (increment) |- (update);
\end{tikzpicture}
\caption{Flow diagram of the Paraexp-EBK method for nonlinear PDEs. The algorithm stops after $K$ iterations, after which the solution is assumed to be converged.}
\label{fig:flowParallel}
\end{figure}

The parallel effiency of the EBK method is achieved by incorporating the parallelization within the iterations that are required to solve nonlinear IVPs. The differences between the serial and the parallel algorithm are illustrated as flow diagrams in Figures~\ref{fig:flowSerial} and~\ref{fig:flowParallel}.

We consider the parallel efficiency in an idealized setting, and estimate a theoretical upper bound here, assuming that communication among the processors can be carried out very efficiently. The computational cost can then be simply estimated by the number of iterations required for the numerical solution to converge. Suppose the computation time of a single EBK iteration is $\tau_0$.
The maximum parallel speedup is then,
\[
\text{speedup} = \frac{K_1 \tau_0}{K_P \tau_0 / P} = \frac{K_1 P}{K_P}
\]
where $K_1$ is number of iterations for serial time integration, and $K_P$ for parallel time integration. The theoretical upper bound of the parallel efficiency is then,
\begin{equation}
\text{efficiency} = \frac{\text{speedup}}{P} = \frac{K_1}{K_P}.
\label{eq:efficiency}
\end{equation}
High parallel efficiency can be achieved if the parallelization does not slow down the convergence of the numerical solution. As will be demonstrated in Section~\ref{sec:burgers}, we typically observe that $K_P$ is not significantly larger than $K_1$, and a near-optimal efficiency is achieved in various relevant cases. For comparison, the parallel efficiency of the Parareal algorithm is formally bounded by $1/K_P$~\cite{Minion:2011aa}. In our case, this upper bound is improved by a factor $K_1$. Parareal is an iterative method for the parallelization of sequential time integration methods, whereas the EBK method, for nonlinear problems, is an iterative method to start with, and its parallelization does not necessarily increase the total number of iterations. Note the PFASST method~\cite{Emmett:2012aa} has a similar estimate of parallel efficiency as the PEBK method.

\begin{figure}[hbtp]
\begin{mdframed}
\textsc{Algorithm parallel time integration.\\}
Given: $A$, $\mbf{u}_0$, $\mbf{g(t,\mbf{u}(t))}$, ...\\
Solve: $\mbf{u}'(t) = A\mbf{u}(t) + \mbf{g}(t,\mbf{u}(t)), \quad \mbf{u}(0) = \mbf{u}_0.$\vspace*{5pt}\\
\textbf{for} $k = 1,\ldots,K$\\
\hspace*{1em}\textbf{for} $j = 1,\ldots,P$ (in parallel)\\
\hspace*{2em}Calculate time-averaged Jacobian matrix $J_k(t)$, see~\eqref{eq:piecewiseJac}.\\
\hspace*{2em}Solve nonhomogeneous part linearized ODE~\eqref{eq:ivpSubproblem}, for $t \in [T_{j-1},T_j]$.\\
\hspace*{2em}Solve homogeneous part linearized ODE~\eqref{eq:ivpSubproblem}, for $t \in [T_j,T_P]$.\\
\hspace*{1em}\textbf{end for}\\
\hspace*{1em}Update solution $\mbf{u}_k(t)$ following~\eqref{eq:updateSolution}.\\
\textbf{end for}
\end{mdframed}
\caption{The algorithm for solving nonlinear problems with the Paraexp exponential block Krylov (PEBK) method, see also Fig.~\ref{fig:flowParallel}.}
\label{fig:algorithm}
\end{figure}

\section{Advection-diffusion equation}
\label{sec:ade}

In this section, we present results of numerical tests where the space-discretized advection-diffusion equation is solved with the EBK method. We demonstrate the consistency and stability of the EBK method for the linear advection-diffusion equation, which is a PDE describing a large variety of transport phenomena~\cite{Hundsdorfer:2013aa}. The spatial discretization, using central finite differencing for sake of illustration, of the PDE yields a linear system of ODEs. The time integration can be parallelized as described in Section~\ref{sec:parlinear}. We illustrate the principle of parallelization for different physical parameters of the advection-diffusion equation, before we move on to nonlinear PDEs in Section~\ref{sec:burgers}. Convergence of the numerical solution is observed for different values of the physical parameters.

In our implementation of the EBK method, the IVP projected onto the block Krylov subspace is solved with the function \texttt{ode15s} in MATLAB (2013b). The relative tolerance of the PEBK solver is denoted by \texttt{tol}. In our tests, the block Krylov subspace is restarted every 20 iterations. For the truncated SVD approximation~\eqref{eq:svd_appr}, $\mbf{p}(t)$ are chosen to be piecewise cubic polynomials, although other types of approximations are possible as well, and are not crucial for the performance of the PEBK method. The sample points per subinterval are Chebyshev nodes. This is not crucial but gives a slightly better approximation than with uniform sample points.

\subsection{Homogeneous PDE}

We consider the advection-diffusion equation with a short pulse initial condition, to clearly distinguish the seperate effects of advection and diffusion. The PDE and the periodic boundary conditions are as follows
\begin{equation}
\begin{aligned}
	u_t + a u_x &= \nu u_{xx}, \quad x \in [0,1], \: t \in [0,1],\\
	u(x,0) &= \sin^{20}\left( \pi x \right),\\
	u(0,t) &= u(1,t),\\
	u_x(0,t) &= u_x(1,t),
\end{aligned}
\label{eq:pulseADE}
\end{equation}
where $a \in \mathbb{R}$ is the advection velocity, and $\nu \in \mathbb{R}$ the diffusivity coefficient. Both parameters are constant in space and time. The PDE is first discretized in space with a second-order central finite difference scheme~\cite{Morton:2005aa}. The corresponding semi-discrete system of ODEs is then
\begin{equation}
\begin{aligned}
	\mbf{v}'(t) &= A\mbf{v}(t) + A\mbf{u}_0,\label{eq:ivpAdvDifEq}\\
	\mbf{v}(0) &= \mbf{0},
	\end{aligned}
\end{equation}
where the matrix $A$ results from the discretization of the spatial derivatives, and represents both the diffusive and the advective term. In~\eqref{eq:ivpAdvDifEq}, the substitution $\mbf{u} = \mbf{v} + \mbf{u}_0$ has been applied, with $\mbf{u}(t)$ being the vector function containing the values of the numerical solution on the mesh at time $t$, with $\mbf{u}_0 = \mbf{u}(0)$. The substitution leads to homogeneous initial conditions. In this case, the source term is constant in time, and its SVD polynomial approximation~\eqref{eq:svd_appr} is exact with $m=1$. This allows us to focus on the two remaining parameters: the grid resolution and the tolerance of the EBK solver.

\begin{figure}[hbtp]
\centering
\begin{subfigure}[t]{0.48\textwidth}
	\includegraphics[width=\linewidth]{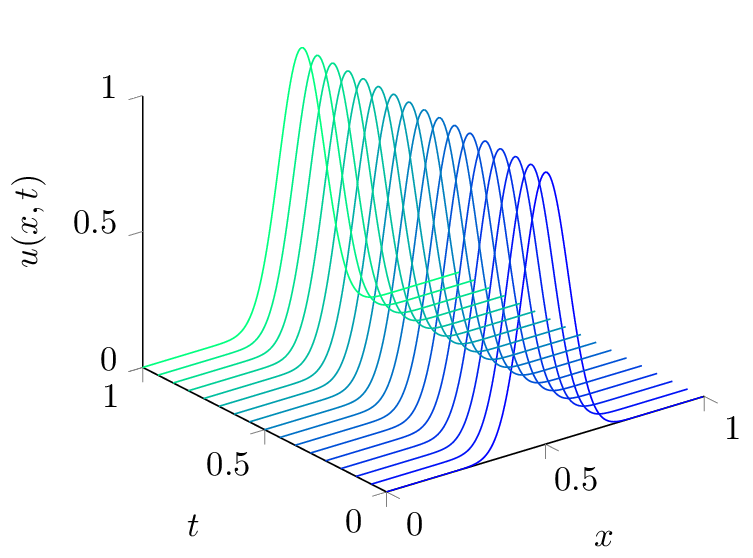}
	\caption{$\mbf{u}_0$}
\end{subfigure}
~
\begin{subfigure}[t]{0.48\textwidth}
	\includegraphics[width=\linewidth]{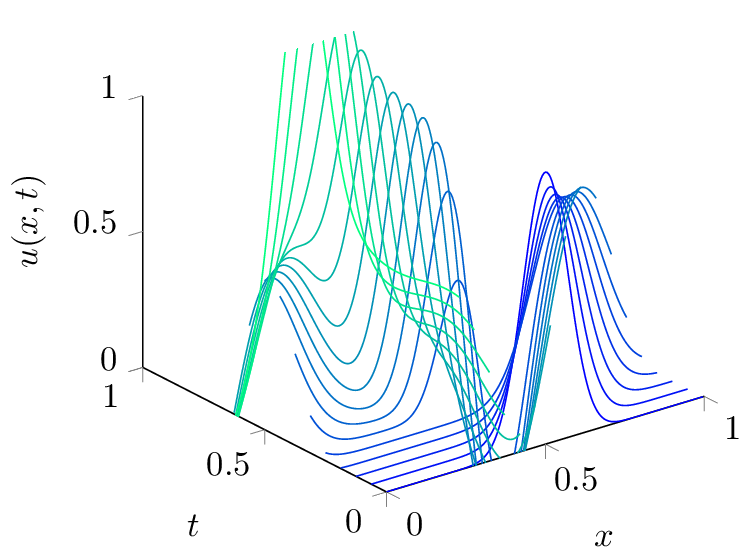}
	\caption{$\mbf{u}_0 + \mbf{v}_1(t)$}
\end{subfigure}
\vspace{5pt}\\
\begin{subfigure}[t]{0.48\textwidth}
	\includegraphics[width=\linewidth]{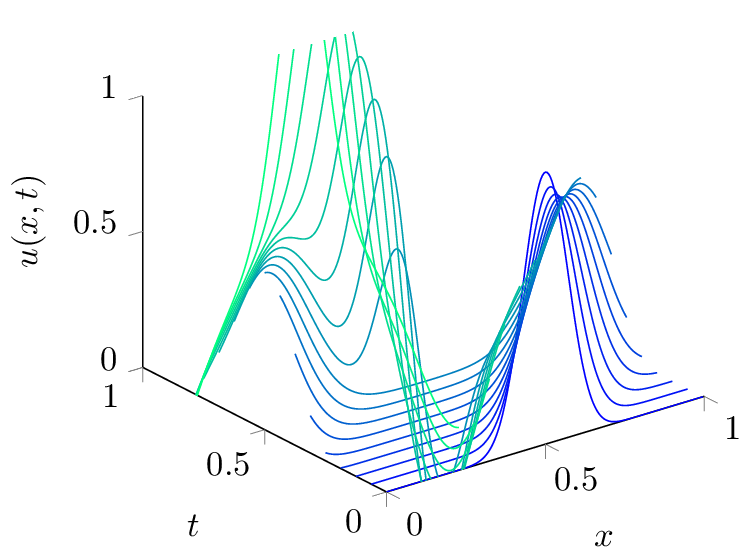}
	\caption{$\mbf{u}_0 + \mbf{v}_1(t) + \mbf{v}_2(t)$}
\end{subfigure}
~
\begin{subfigure}[t]{0.48\textwidth}
	\includegraphics[width=\linewidth]{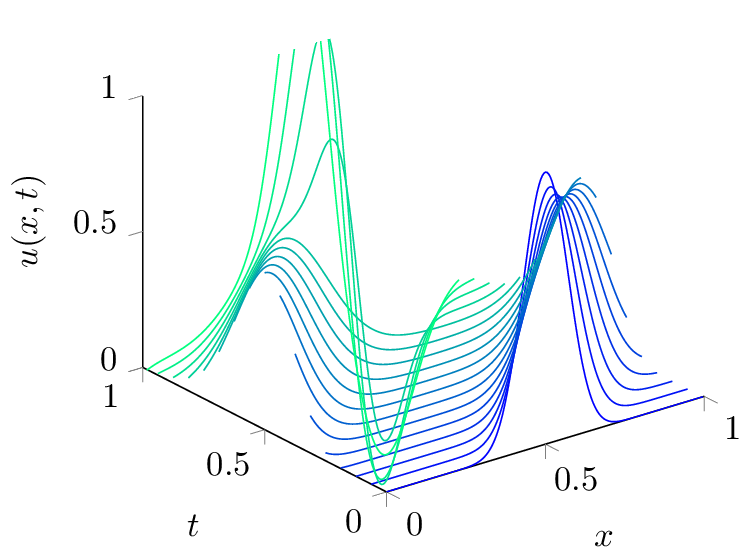}
	\caption{$\mbf{u}_0 + \mbf{v}_1(t) + \mbf{v}_2(t) + \mbf{v}_3(t)$}
\end{subfigure}
\vspace{5pt}\\
\begin{subfigure}[t]{0.48\textwidth}
	\includegraphics[width=\linewidth]{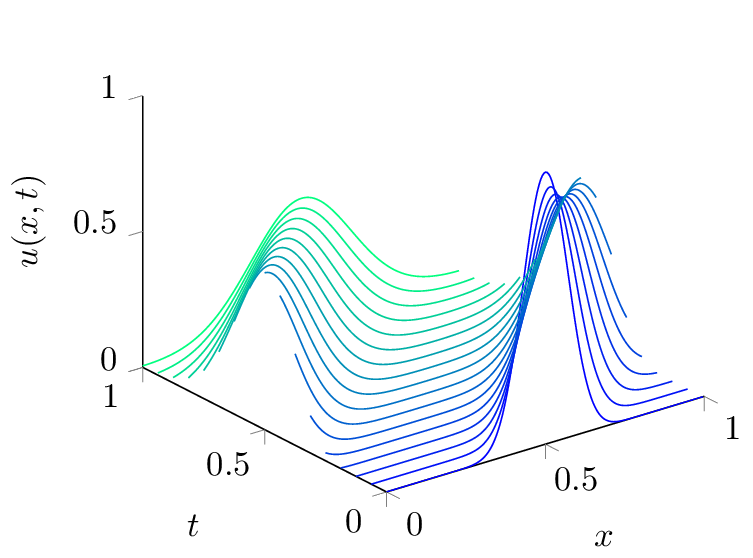}
	\caption{$\mbf{u}(t) = \mbf{u}_0 + \sum_{j=1}^4\mbf{v}_j(t)$}
\end{subfigure}

\caption{Superposition of subsolutions to the advection-diffusion equation, with partitioning $\mathbf{T} = \{0,0.25,0.5,0.75,1\}$, $\Delta x = 5 \cdot 10^{-3}$, $a = 1$, and $\nu = 10^{-2}$. The color depends on the time coordinate: blue corresponds to $t=0$ and green to $t=1$.}
\label{fig:superposition}
\end{figure}

The linear IVP~\eqref{eq:ivpAdvDifEq} is decoupled into independent subproblems by partitioning of the source term ($A\mbf{u}_0$) on the time interval. The superposition of the subsolutions is illustrated in Fig.~\ref{fig:superposition}, in which the time interval of interest, $[0,1]$, has been partitioned into four equal subintervals, to be integrated on four processors. The sum of the initial condition and the subsolutions gives the final solution on the entire time interval, see Section~\ref{sec:parlinear}.

\begin{figure}[hbtp]
\centering
\begin{subfigure}[t]{\textwidth}
	\includegraphics[scale=1]{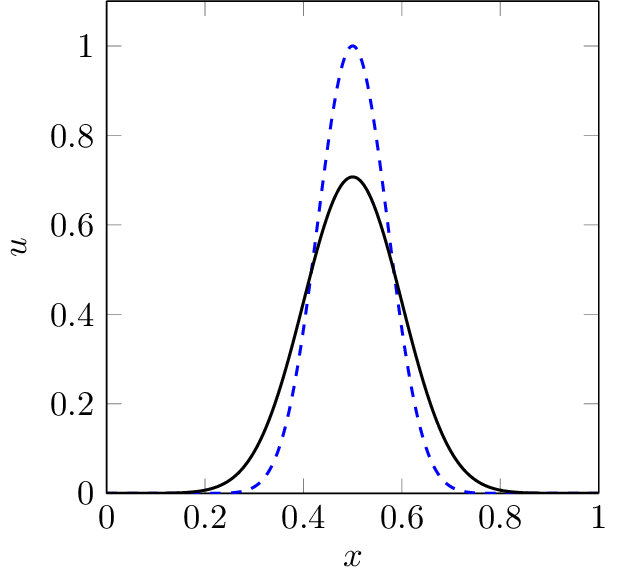}~
	\includegraphics[scale=1]{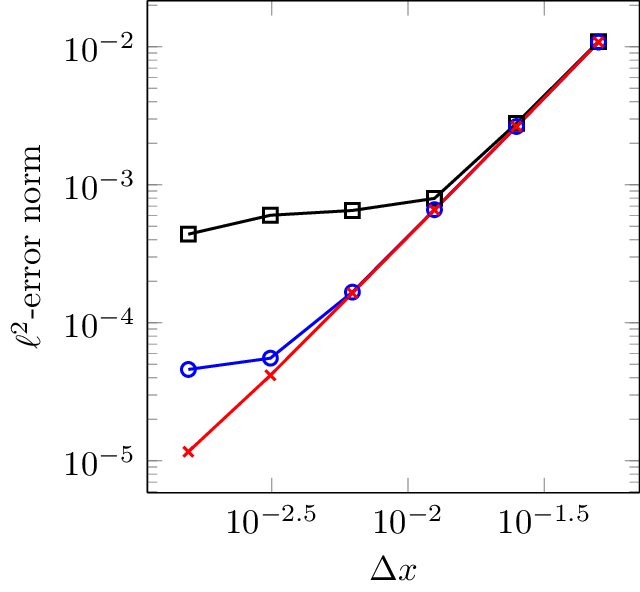}
	\caption{Velocity $a=0$ and diffusivity coefficient $\nu = 10^{-2}$.}
	\label{fig:solutionDif}
\end{subfigure}

\begin{subfigure}[t]{\textwidth}
	\includegraphics[scale=1]{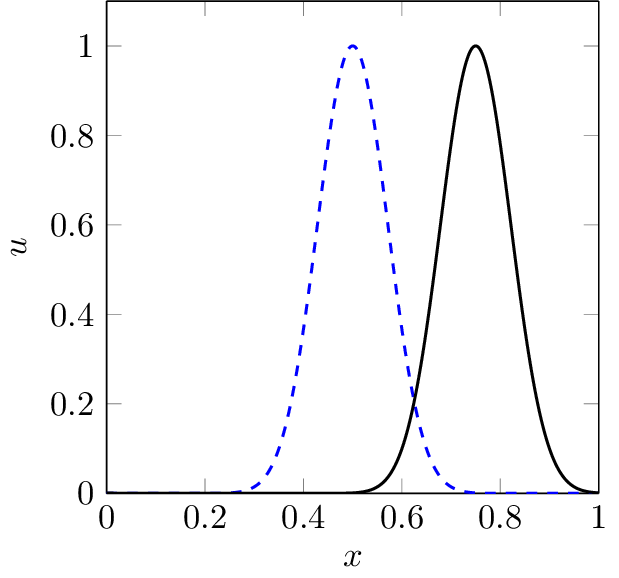}~
	\includegraphics[scale=1]{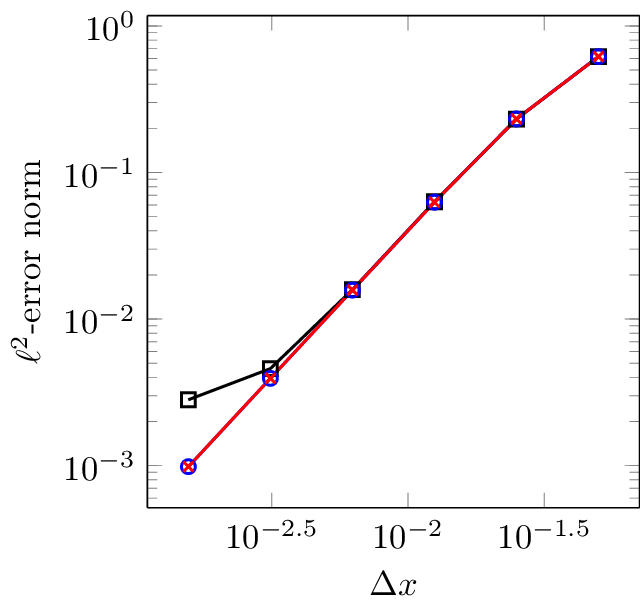}
	\caption{Velocity $a=1$ and diffusivity coefficient $\nu = 0$.}
	\label{fig:solutionAdv}
\end{subfigure}

\begin{subfigure}[t]{\textwidth}
	\includegraphics[scale=1]{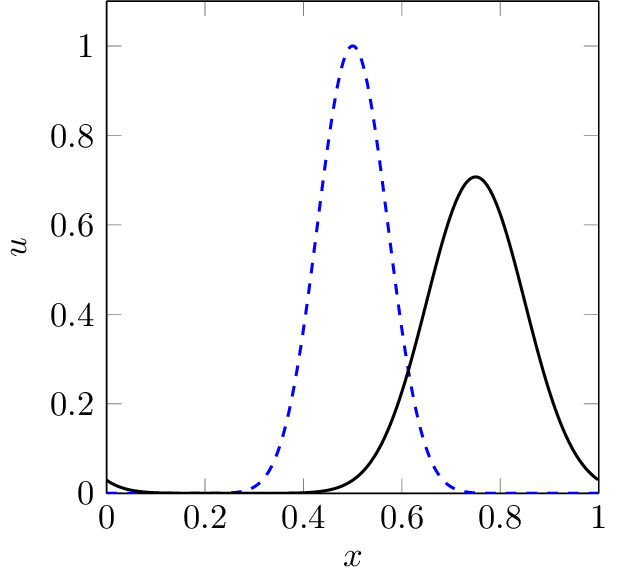}~
	\includegraphics[scale=1]{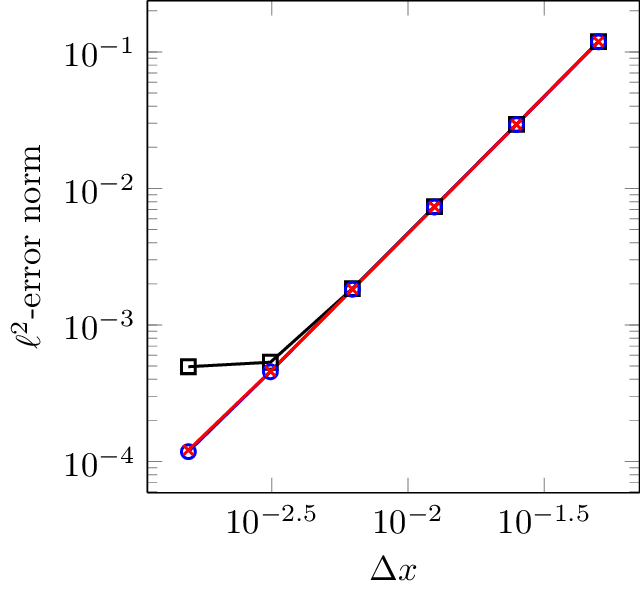}
	\caption{Velocity $a=1$ and diffusivity coefficient $\nu = 10^{-2}$.}
	\label{fig:solutionAde}
\end{subfigure}

\caption{Left: numerical solution on the finest mesh at $t = 0$ (dashed) and $t=0.25$ (solid). Note that the solution for $a=1$ and $\nu = 0$ does not show odd-even decoupling because of the high spatial resolution used. Right: convergence of the numerical solution (at final time $t=1$), for $\mathtt{tol}$: $\square\:10^{-2}$, ${\color{blue} \circ}\:10^{-4}$, ${\color{red} \times}\:10^{-6}$.}
\label{fig:solnADE}
\end{figure}

In the following numerical experiments, $P = 8$ subintervals are used. The advection-diffusion equation is solved for three different combinations of $a$ and $\nu$, see Fig.~\ref{fig:solnADE}. The solutions are computed for different $\mathtt{tol}$ and $\Delta x$. The discretization error is controlled by the time integration with $\mathtt{tol}$, and by the spatial discretization with $\Delta x$. The error of the numerical solution is measured in the relative $\ell^2$-norm,
\begin{equation}
\Vert \mbf{u}_h(T) - \mbf{u}(T)\Vert / \Vert \mbf{u}(T) \Vert,
\label{eq:relerr}
\end{equation}
where the exact solution $\mbf{u}(T)$ is on the mesh, i.e., it has the entries $u_j(T) = u(x_j,T)$. The exact solution of Eq.~\eqref{eq:pulseADE} is
\begin{equation}
u(x,t) = \sum^{20}_{n=0} a_n \cos \left(n \pi (x - a t)\right) \exp\left(-(n\pi)^2\nu t\right),
\end{equation}
with coefficients
\begin{align}
a_0 &= \frac{1}{2} \int_{-1}^1 \sin^{20} (\pi x) \, dx,\\
a_n &= \int_{-1}^1 \sin^{20} (\pi x) \cos (n \pi x) \, dx .
\end{align}
The coefficients are calculated using numerical quadrature with high precision. The error in~\eqref{eq:relerr} shows second-order convergence with $\Delta x$, as expected from the truncation error of the finite difference scheme. The convergence plots show that we are able to control the error of the parallel time integration method. In this case, the final error can be made to depend only on the spatial resolution and the tolerance of the EBK method. Also, the EBK method has no principal restrictions on the timestep size, as it directly approximates the exact solution of Eq.~\eqref{eq:ivpAdvDifEq} by Krylov subspace projections.

\clearpage

\subsection{Parallel efficiency}
\label{sec:ade_efficiency}

\begin{figure}
\centering
\begin{minipage}{.48\textwidth}
\centering
\includegraphics[width=\linewidth]{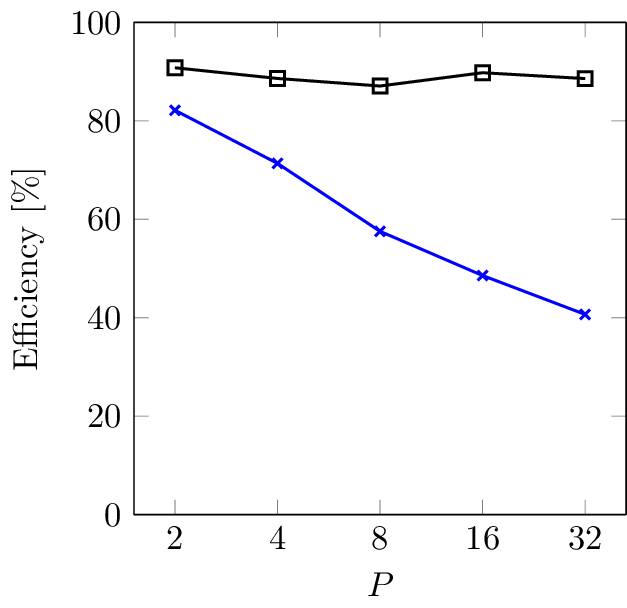}
\caption{Parallel efficiency for solving the advection-diffusion equation. $\square$~PEBK; $\color{blue}\times$~Paraexp/CN.}
\label{fig:parallel_efficiency}
\end{minipage}%
\hfill
\begin{minipage}{.48\textwidth}
\centering
\includegraphics[width=\linewidth]{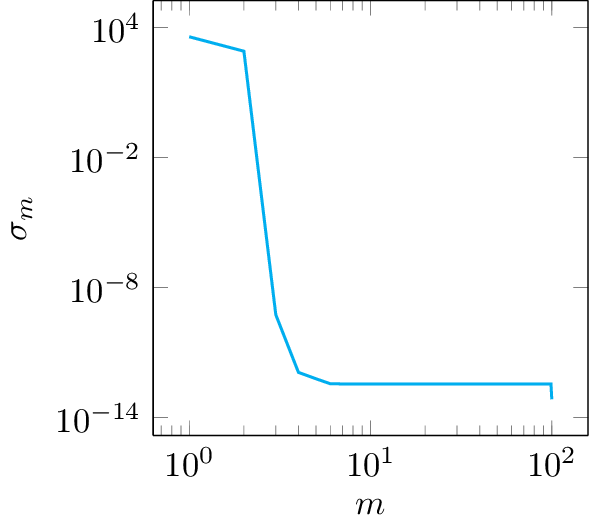}
\caption{Singular values of the matrix composed of the source term samples.}
\label{fig:singular_values}
\end{minipage}
\end{figure}

In the previous examples, we solved homogeneous IVPs. Nonhomogeneous problems are generally more expensive to solve, because an accurate SVD approximation of the source term~\eqref{eq:svd_appr} requires more singular values, which increases the block width of the block Krylov subspace~\eqref{eq:blockKrylovSubspace}. Parallel speedup can then be expected by splitting the nonhomogeneous problem into subproblems~\eqref{eq:linivpSubproblem}, which require less individual effort to solve by the PEBK method. To measure the parallel efficiency of our algorithm for nonhomogeneous IVPs, we introduce a source term $f(x,t)$ in the advection-diffusion equation
\begin{equation}
u_t + a u_x = \nu u_{xx} + f(x,t),
\label{eq:ade}
\end{equation}
with $\nu = 10^{-2}$, and $a = 1$. The source term is chosen such that the solution is a series of five travelling pulses
\begin{equation}
u(x,t) = \tfrac{1}{2} - \tfrac{1}{2} \cos\left( 10 \pi (x - t) \right).
\label{eq:manusol}
\end{equation}
The mesh width of the spatial discretization is $\Delta x = 10^{-3}$, such that the error due to the spatial discretization is small compared to the time integration error. The mesh width gives a semi-discrete system with a $1000 \times 1000$ matrix, which is a suitable problem size for testing the EBK method. The tolerance of the EBK method is set to $10^{-4}$. The SVD polynomial approximation is constructed from 100 time samples per subinterval. {The singular values are plotted in Fig.~\ref{fig:singular_values}, which reveals that the first two singular values are several orders larger than the rest. Therefore, we retain only the first two singular values in the truncated SVD. The decay of singular values is guaranteed by the upper bound from Theorem~\ref{The1}.} We have verified that the SVD approximation is sufficiently accurate, i.e., the approximation error, measured in the $L^2$-norm, is less than the tolerance of the EBK method.

We compare the parallel efficiency of the Paraexp-EBK implementation with a convential implementation of the Paraexp algorithm~\cite{Gander:2013aa}, where the nonhomogeneous problem is integrated with the Crank--Nicolson (CN) scheme. The linear system is solved directly using MATLAB. The matrix exponential propagator, for the homogeneous problem in the Paraexp method, is realized with an Arnoldi method using the shift-and-invert technique~\cite{Eshof:2006aa}.

In order to have a Courant number of one, we take $\Delta t = 10^{-3}$ (for $\Delta x = 10^{-3}$ and $a=1$). According to the Paraexp method, the time step size needs to be decreased in parallel computation by a factor ${P}^{1/2q}$~\cite{Gander:2013aa}, where $q$ is the order of the time integration method, in order to control the error. In case of the Crank--Nicolson scheme, we have $q = 2$.

As discussed in Section~\ref{sec:parlinear}, there is no communication required between processors, except for the superposition of the solutions to the subproblems at the end. Therefore, we are able to test the parallel algorithm on a serial computer by measuring the computation time of each independent subproblem. The computation time of the serial time integration is denoted $\tau_0$. For the parallel time integration, we measure the computation times of the nonhomogeneous and the homogeneous part of the subproblems separately, denoted by $\tau_1$ and $\tau_2$ respectively. The parallel speedup can then be estimated as
\begin{equation}
\text{speedup} = \frac{\tau_0}{\max(\tau_1) + \max(\tau_2)},
\end{equation}
where we take the maximum value of $\tau_1$ and $\tau_2$ over all parallel processes. The timings of the EBK method and Paraexp are listed in Table~\ref{tab:efficiency}.

The parallel efficiency of both methods is illustrated in Figure~\ref{fig:parallel_efficiency}. Note that the parallel efficiency of the standard Paraexp method steadily decreases, whereas the PEBK method maintains a constant efficiency level around $90\%$. The decrease in efficiency of the standard Paraexp method is due to the reduced time step size in the Crank--Nicolson scheme, with respect to its serial implementation. The nonhomogeneous part of the subproblems requires more computation time as the number of processors increases.

The numerical test confirms the initial assumption that parallel speedup with the EBK method can be achieved by decomposing the originial problem into simpler independent subproblems~\eqref{eq:linivpSubproblem}. Furthermore, the Paraexp implementation using the EBK method, appears to be more efficient than one using a traditional time-stepping method.


\begin{table}[hbtp]
\caption{Parallel effiency of the Paraxp-EBK method and the Paraexp/Crank--Nicolson method for the advection-diffusion equation, with number of processors $P$. Timing $\tau_0$ corresponds to the serial algorithm. For the parallel algorithm, $\tau_1$ and $\tau_2$ are timings of the nonhomogeneous and homogeneous problem respectively.}
\centering
\begin{tabular}{l c c c c c c c}
\toprule
& & \multicolumn{2}{c}{Serial} & \multicolumn{3}{c}{Parallel} & \\
\cmidrule(lr){3-4}
\cmidrule(lr){5-7}
& $P$ & $\tau_0$ & Error & $\max{(\tau_1)}$ & $\max{(\tau_2)}$ & Error & \small \emph{Efficiency}\\
\midrule
PEBK & 2 & 1.43e+00 & 3.05e-04 & 7.22e-01 & 6.44e-02 & 2.70e-04 & 91 $\%$ \\ 
 & 4 & 2.81e+00 & 3.05e-04 & 7.23e-01 & 7.06e-02 & 2.68e-04 & 89 $\%$ \\ 
 & 8 & 5.59e+00 & 3.05e-04 & 7.40e-01 & 6.21e-02 & 2.88e-04 & 87 $\%$ \\ 
 & 16 & 1.14e+01 & 3.05e-04 & 7.30e-01 & 6.28e-02 & 3.22e-04 & 90 $\%$ \\ 
 & 32 & 2.27e+01 & 3.05e-04 & 7.33e-01 & 6.63e-02 & 3.65e-04 & 89 $\%$ \\ 
\midrule
Paraexp & 2 & 2.15e+00 & 4.56e-04 & 1.29e+00 & 1.66e-02 & 4.10e-04 & 82 $\%$ \\ 
 & 4 & 4.37e+00 & 4.56e-04 & 1.52e+00 & 1.28e-02 & 3.80e-04 & 71 $\%$ \\ 
 & 8 & 8.56e+00 & 4.56e-04 & 1.85e+00 & 1.24e-02 & 3.59e-04 & 58 $\%$ \\ 
 & 16 & 1.70e+01 & 4.56e-04 & 2.18e+00 & 1.31e-02 & 3.43e-04 & 49 $\%$ \\ 
 & 32 & 3.43e+01 & 4.56e-04 & 2.62e+00 & 1.16e-02 & 3.32e-04 & 41 $\%$ \\ 
\bottomrule
\end{tabular}
\label{tab:efficiency}
\end{table}

\clearpage

\section{Burgers equation}
\label{sec:burgers}

In the previous section, the PEBK method was applied to a linear PDE. In this section, the performance of the PEBK method is tested on a nonlinear PDE, the viscous Burgers equation.

\subsection{Travelling wave}
\label{sec:burgersWave}
Consider the viscous Burgers equation,
\begin{align}
&u_t + u u_x = \nu u_{xx} + f(x,t), \quad x \in [0,1], t \in [0,1], \label{eq:burgersEq}
\end{align}
where $\nu \in \mathbb{R}$ denotes the diffusivity (or viscosity) coefficient. In the following experiments, we take $\nu = 10^{-2}$, such that the problem is dominated by the nonlinear convective term. As in the previous example, the boundary conditions are periodic. Exact solutions to the Burgers equation can be found by the Cole--Hopf transformation~\cite{Hopf:1950aa}. In this test case, we construct a desired solution by introducing an appropriate source term, as shown for the advection-diffusion equation in Section~\ref{sec:ade_efficiency}. The the source term balances the dissipation of energy due to diffusion, and prevents the solution of vanishing in the limit $t \rightarrow \infty$.

\begin{figure}[hbtp]
\centering
\begin{subfigure}[t]{0.48\textwidth}
	\includegraphics[scale=1]{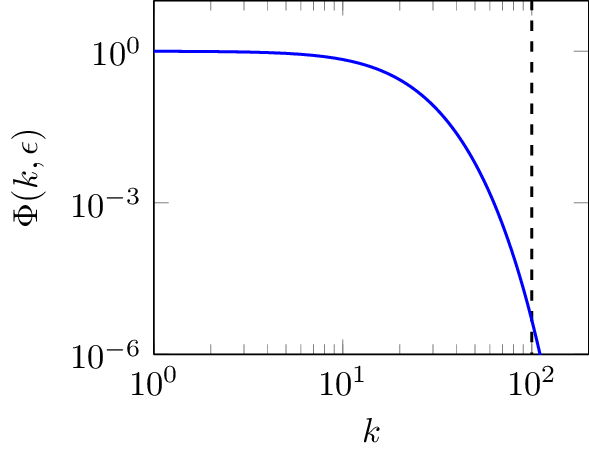}
	\caption{Smoothing function (solid line) and cutoff wavenumber (dashed line).}
	\label{fig:smoothingFunction}
\end{subfigure}
~
\begin{subfigure}[t]{0.48\textwidth}
	\includegraphics[scale=1]{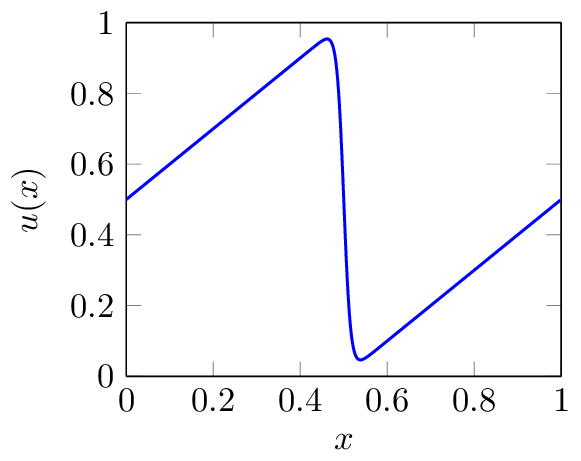}
	\caption{Solution for $\epsilon = 0.1$ and $k_{max}=100$.}
	\label{fig:sawtooth}
\end{subfigure}
\vspace{-.5cm}
\caption{Manufactured solution~\eqref{eq:modifiedSawtooth} of the Burgers equation.}
\end{figure}

The Burgers equation is an important equation in fluid dynamics. It can be seen as a simplification of the Navier--Stokes equation, which retains the nonlinear convective term, $u u_x$. In the limit $\nu \rightarrow 0$, the nonlinearity may produce discontinuous solutions, or shocks, so that a typical solution may resemble a sawtooth wave. A sawtooth wave can be represented by an infinite Fourier series. A smooth version of the sawtooth wave can be obtained by the modified Fourier series,
\begin{equation}
u(\xi) = \frac{1}{2} - \sum_{k = 1}^{k_{\text{max}}} \frac{1}{\pi k} \sin( 2 \pi k \xi ) {\Phi(k,\epsilon)} , \quad \xi \in \mathbb{R}, \label{eq:modifiedSawtooth}
\end{equation}
where $k$ is the wavenumber, $k_{\text{max}}$ a cutoff wavenumber, and $\Phi(k,\epsilon)$ is a smoothing function, which supresses the amplitudes associated to high wavenumbers:
\begin{equation}
\Phi(k,\epsilon) = \left[ \frac{\pi k \epsilon /2}{\sinh(\pi k \epsilon /2)}\right],
\end{equation}
Here, $\epsilon$ is the smoothing parameter. The smoothing function $\Phi(k,\epsilon)$ is motivated by the viscosity-dependent inertial spectrum of the Burgers equation found by Chorin \& Hald~\cite{Chorin:2005aa}. The smoothing function for $\epsilon = 0.1$ is shown in Figure~\ref{fig:smoothingFunction}. As the smoothing function rapidly decreases with wavenumber, we choose a cutoff wavenumber of $k_{\text{max}} = 100$. The value $\epsilon = 0.1$ is found to produce a smooth version of a sawtooth wave, as shown in Figure~\ref{fig:sawtooth}. 

We consider a wave travelling in the positive $x$-direction by introducing the  parametrization $\xi = x - \frac{1}{2}t + \frac{1}{2}$ in~\eqref{eq:burgersEq}. The source term is then readily found by substituting the chosen solution~\eqref{eq:modifiedSawtooth} into the Burgers equation~\eqref{eq:burgersEq},
\begin{equation}
f(x,t) = u_t + u u_x - \nu u_{xx}.
\end{equation}
The space derivatives in the Burgers equation are discretized with a second-order central finite difference schemes, using a uniform mesh width $\Delta x > 0$. The error of the numerical solution is again measured in the relative $\ell^2$-norm, cf.~\eqref{eq:relerr}.

The tolerance of the PEBK method is set to $10^{-4}$. The SVD approximation of the source term, which includes the nonlinear term, is constructed from $s=50$ samples per subinterval, and reveals that $m=12$ singular values are sufficient in the truncated SVD, so that the error of the SVD approximation is less than the tolerance of the PEBK method.

\begin{figure}[hbtp]
\centering
\includegraphics[scale=1]{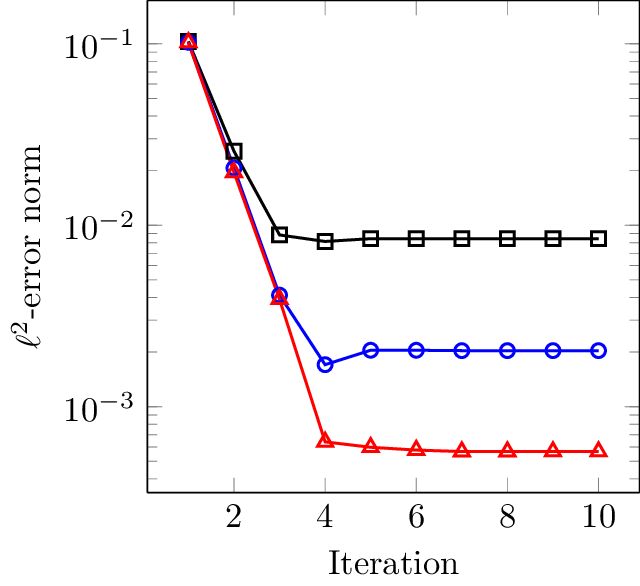}
\caption{Error history with $\Delta T = 0.1$. $\square$, $\Delta x = 10^{-2}$; $\color{blue} \circ$, $\Delta x = 5 \cdot 10^{-3}$; $\color{red} \triangle$, $\Delta x = 2.5 \cdot 10^{-3}$.}
\label{fig:error_differentdx}
\end{figure}

The nonlinear system of ODEs is solved iteratively, as outlined in Section~\ref{sec:ebk_nonlinear}. Figure~\ref{fig:error_differentdx} shows the error history at different grid resolutions. Here, the error converges to a value that depends on the mesh width. In other words, the final error of the time integration method is much smaller than the error due to the spatial discretization.

\begin{figure}[hbtp]
\centering
\begin{subfigure}[t]{0.48\textwidth}
	\includegraphics[scale=1]{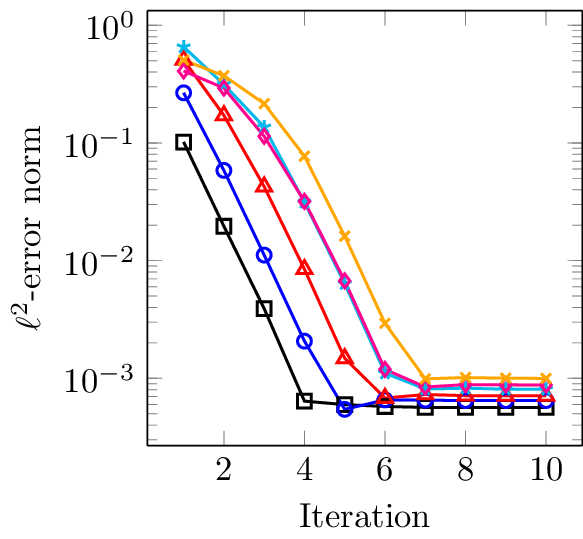}
	\caption{$\Delta T = 0.1$, $T = P\Delta T$.}
	\label{fig:error_dtfix}
\end{subfigure}
~
\begin{subfigure}[t]{0.48\textwidth}
	\includegraphics[scale=1]{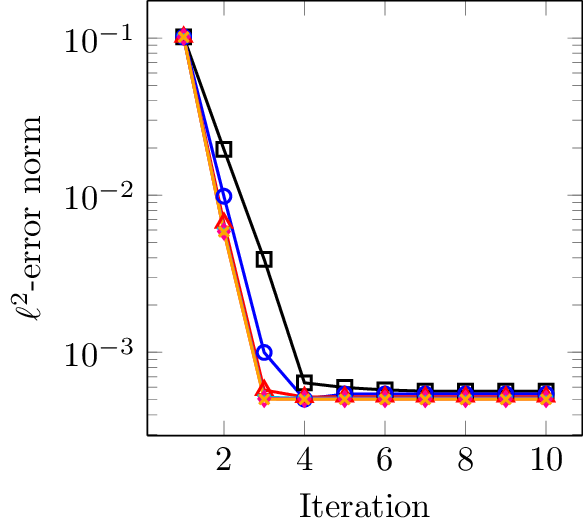}
	\caption{$\Delta T = T/P$, $T = 0.1$.}
	\label{fig:error_tendfix}
\end{subfigure}
\vspace{-.5cm}
\caption{Error history for the Paraxp-EBK method for a fixed subinterval size and a fixed final time. $\square$, $P=1$; $\color{blue} \circ$, $P=2$; $\color{red} \triangle$, $P=4$; $\color{cyan} \star$, $P=8$; $\color{magenta} \diamond$, $P=16$; $\color{orange} \times$, $P=32$.}
\label{fig:be1_history}
\end{figure}

The PEBK method is parallelized as discussed in Section~\ref{sec:parnonlinear}. The time interval $[0,T]$ can be partitioned into $P$ subintervals with a uniform subinterval size $\Delta T$. In the following experiment, the complete time interval is extended with each subinterval added, $T = P\Delta T$. The mesh width is $\Delta x = 2.5 \cdot 10^{-3}$. Figure~\ref{fig:error_dtfix} shows that increasing the number of processors and hence the simulated time, generally increases the number of iterations required to achieve the same level of accuracy. The increased number of required iterations implies a decrease in theoretical parallel efficiency, which can be estimated by the ratio $K_1/K_P$, see~\eqref{eq:efficiency}.

Next, the final time is kept constant, and the subinterval size reduces with an increase in the number of processors, $\Delta T = T/P$. Figure~\ref{fig:error_tendfix} shows that the number of iterations is roughly independent of the number of processors, i.e., in this case, the parallel efficiency, $K_1/K_P$, does not decrease with $P$. Parallel speedup might also be improved by the fact that the SVD approximation converges faster on smaller subintervals, see Theorem~\ref{The1}. That is, fewer singular values have to be retained in the truncated SVD in order to achieve a certain accuracy. Also, the number of samples of the source term could be decreased on smaller subintervals.

\subsection{Parallel efficiency}

\begin{figure}[hbtp]
\centering
\includegraphics[scale=1]{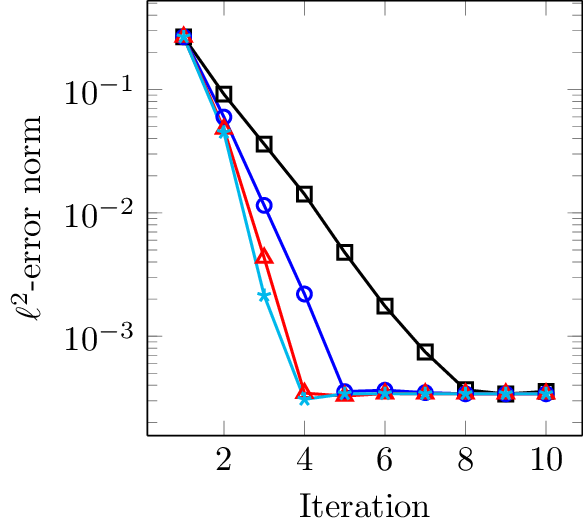}
\caption{Error history for the Paraxp-EBK method with $\Delta T = 0.2/P$. $\square$, $P=1$; $\color{blue} \circ$, $P=2$; $\color{red} \triangle$, $P=4$; $\color{cyan} \star$, $P=8$.}
\label{fig:be_history}
\end{figure}

\begin{figure}[hbtp]
\centering
\begin{subfigure}[t]{0.48\textwidth}
	\includegraphics[width=\linewidth]{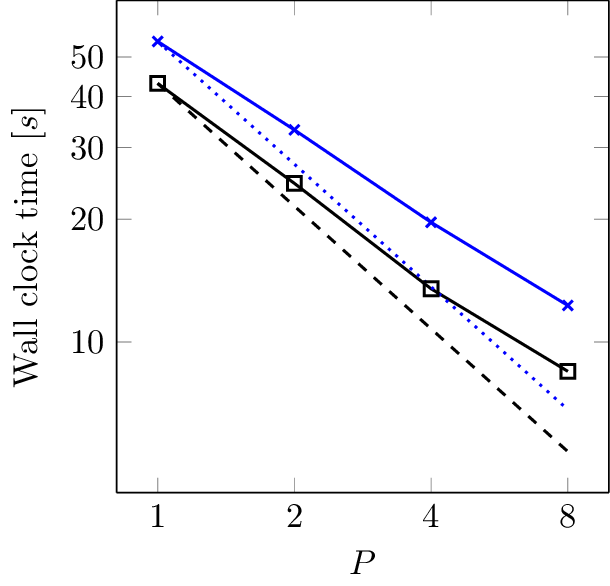}
	\caption{$\color{blue} \times$: $\nu = 10^{-1}$; $\square$: $\nu = 10^{-2}$; {\color{blue} dotted line}: $\nu=10^{-1}$ (ideal); dashed line: $\nu=10^{-2}$ (ideal).}
	\label{fig:be_timing}
\end{subfigure}
~
\begin{subfigure}[t]{0.48\textwidth}
	\includegraphics[width=\linewidth]{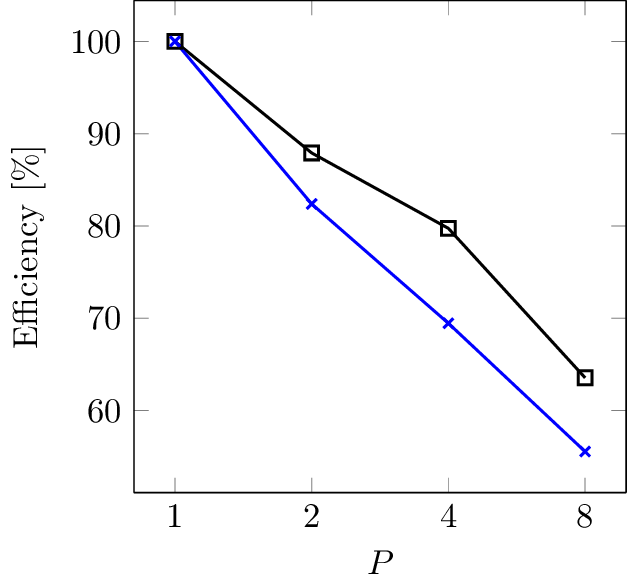}
	\caption{$\color{blue} \times$: $\nu = 10^{-1}$; $\square$: $\nu = 10^{-2}$.}
	\label{fig:be_efficiency}
\end{subfigure}
\vspace{-.5cm}
\caption{Total computation time (left) and parallel efficiency (right) of ten PEBK iterations, with $\Delta T = 0.2/P$.}
\end{figure}

{
In this section, we test the parallel efficiency of the PEBK method for the viscous Burgers equation. In the following experiments, the source term is chosen such that the solution of the spatially discretized system is equal to the exact solution described in Section~\ref{sec:burgersWave}. In other words, the source term accounts for the error by the spatial discretization, and we measure only the error due to the time integration. The spatial discretization is here performed with a mesh width of $\Delta x = 2 \cdot 10^{-3}$.

If the subinterval size, $\Delta T$, is fixed, the number of iterations generally increases in case the number of processors $P$ increases, see Fig.~\ref{fig:be1_history}. Increasing the number of iterations would reduce the parallel efficiency, according to~\eqref{eq:efficiency}. We therefore fix the final time $T$, such that the subinterval size decreases with increasing $P$. Also, the total number of samples over $[0,T]$ is fixed, such that the local number of samples decreases with increasing $P$. To keep the global distribution of samples constant in the parallel computations, the local sample points are uniformly distributed instead of Chebyshev points as used in previous experiments. In our experiments, we use $\Delta T = 0.2/P$ and $s = 128/P$. Furthermore, the number of (retained) singular values is $m = 12$, and the tolerance of the EBK method is set to $10^{-4}$. The error history for different $P$ is shown in Fig.~\ref{fig:be_history}. The results show that the convergence rate of the PEBK method can improve by increasing $P$. The nonlinear corrections, see \eqref{eq:ivp4}, appear to be more effective on smaller subintervals.

As opposed to the linear problem in Section~\ref{sec:ade_efficiency}, communication between the parallel processes is here required because of the waveform relaxation method. The parallel communication overhead is assumed to be very small, such that the parallel computations can be emulated on a serial computer. Per iteration of the PEBK method, the computation time of each individual process is measured, after which the maximum is stored, i.e., the computation time of the slowest process. The total computation time of the PEBK method is then taken as the sum, over the total number of iterations performed, of the maxima. In this experiment, we have measured a total of ten PEBK iterations.

The computation times for $\nu = 10^{-1}$ and $\nu = 10^{-2}$ are shown in Fig.~\ref{fig:be_timing}, which clearly illustrate a parallel speedup. The slightly higher timings of $\nu = 10^{-1}$ could be attributed to the increased stiffness of the problem. Figure~\ref{fig:be_efficiency} shows that the parallel efficiency of the PEBK method steadily decreases with the number of (virtual) processors, $P$. There is no significant difference between the performance of the method at these two values of the viscosity coefficient. The parallel efficiency might even be further tuned in practice, based on the previous observation that the required number of iterations could decrease with higher $P$. Also, the number of singular values, $m$, could possibly be reduced on smaller subintervals, based on~\eqref{The1}, which would further enhance the potential parallel speedup.
}

The PFASST algorithm~\cite{Emmett:2012aa} shows a good parallel efficiency for the viscous Burgers equation as well. It is unknown whether changing the viscosity coefficient has an impact on the performance of PFASST. Our observations point in the direction of a good parallel efficiency of the PEBK method for simulations of the Navier--Stokes equation at high Reynolds numbers. For these problems there is a high demand for highly efficient parallel solvers~\cite{Poel:2015aa}. The Parareal algorithm was for example reported to perform poorly in advection-dominated problems~\cite{Fischer:2005aa,Steiner:2015aa}.

\subsection{Multiscale example}

\begin{figure}[hbtp]
\centering
\begin{subfigure}[t]{0.48\textwidth}
	\includegraphics[width=\textwidth]{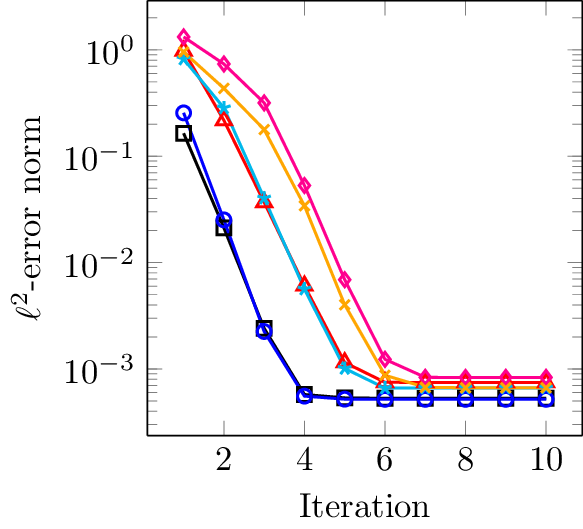}
	\caption{$k_0=4$}
	\label{fig:burgers_k4}
\end{subfigure}
~
\begin{subfigure}[t]{0.48\textwidth}
	\includegraphics[width=\textwidth]{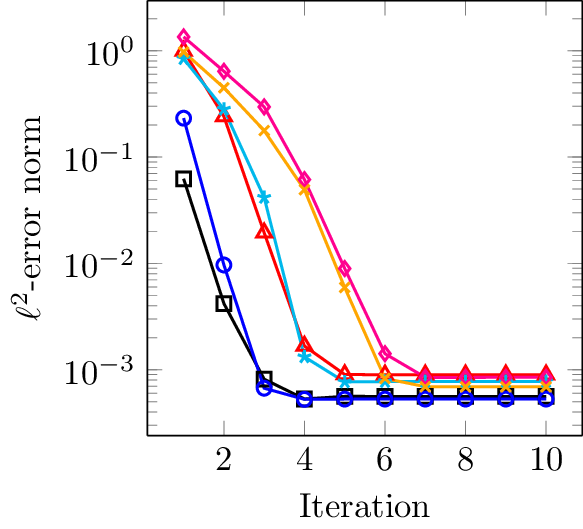}
	\caption{$k_0=16$}
	\label{fig:burgers_k16}
\end{subfigure}
\caption{Error history with $\Delta T = 0.1$ and $T = P \Delta T$. $\square$, $P=1$; $\color{blue} \circ$, $P=2$; $\color{red} \triangle$, $P=4$; $\color{cyan} \star$, $P=8$; $\color{magenta} \diamond$, $P=16$; $\color{orange} \times$, $P=32$.}\label{fig:burgers_sine}
\end{figure}

Turbulent flow is a multiscale phenomenon that is characterized by a wide range of length and time scales. In case of the Burgers equation, we can emulate such a solution by combining two functions that are periodic in both space and time, e.g.,
\begin{align}
u(x,t) = \sin( 2 \pi x) \sin( 2 \pi t) + \frac{1}{k_0}\sin( 2 k_0 \pi x) \sin( 2 k_0 \pi t), \quad k_0 > 1 .
\label{eq:multiscale_solution}
\end{align}
The solution features a large scale mode with wavenumber one, and a smaller scale mode with wavenumber $k_0>1$. This combination allows the construction of an arbitrarily wide dynamic range. The factor $1/k_0$ is included in compliance with an assumed energy distribution of $\left[\hat{u}(k)\right]^2 \propto k^{-2}$~\cite{Chorin:2005aa}, where $\hat{u}(k)$ is the Fourier transform, 
\begin{equation}
\hat{u}(k) = \int_{-\infty}^\infty u(x,t)  \exp( -2\pi i k (x + t) )\,dx\,dt,
\end{equation}
where we have used the fact that the solution is periodic in space and time. The previous experiment, see Fig.~\ref{fig:error_dtfix}, is repeated for the manufactured solution~\eqref{eq:multiscale_solution} with different values of $k_0$. Figure~\ref{fig:burgers_sine} illustrates that widening the spectrum does not significantly affect the convergence of the PEBK iterations. Remarkably, the curves appear to form pairs based on $P$.

\section{Conclusions}
\label{sec:conclusions}
We propose an implementation of the Paraexp method with enhanced parallelism based on the exponential block Krylov (EBK) method. Furthermore, the method, Paraexp-EBK (PEBK), is extended to solve nonlinear PDEs iteratively by a waveform relaxation method. The nonlinear terms are represented by a source term in a nonhomogeneous linear system of ODEs. Each iteration the source term is updated with the latest solution. The convergence of the iterative process can be accelerated by adding a correction term based on the Jacobian matrix of the nonlinear term. Each iteration the initial value problem can then be decoupled into independent subproblems, which can be solved parallel in time. Essentially, we implement the Paraexp method within a waveform relaxation approach in order to integrate nonlinear PDEs. Also, the Paraxp-EBK (PEBK) method is used to integrate both the homogeneous and the nonhomogeneous parts of the subproblems. This is in contrast to the original Paraexp method, which assumes a convential time integration method for the nonhomogeneous parts.

The PEBK method is tested on the advection-diffusion equation for which we demonstrate the parallelization concept for linear PDEs. The parallelization also works in cases without diffusion present, in which the PDE is purely hyperbolic. The parallel efficiency is compared with a Crank--Nicolson (CN) scheme parallelized with the Paraexp algorithm. The parallel efficiency of the PEBK method remains roughly constant around $90\%$. On the other hand, the parallel efficiency of the CN/Paraexp combination steadily decreases with the number of processors.

As a model nonlinear PDE, the viscous Burgers equation is solved. The number of waveform relaxation iterations required for a certain error tolerance increases, when the relative importance of nonlinearity grows by decreasing the viscosity coefficient. Good parallel efficiency of the EBK method was observed for different values of the viscosity coefficient.
Since the nonlinear convective term in the Burgers equation is shared by the Navier--Stokes equation, the presented results give a hint of the potential of the PEBK method as an efficient parallel solver in turbulent fluid dynamics, where nonlinearity plays a key role. The question remains how to treat the pressure in the incompressible Navier--Stokes equation, which enforces the incompressibility constraint on the velocity field (see for possible approaches in~\cite{Edwards:1994aa,Newman:2004aa}). This will be explored in future work.

\bibliographystyle{abbrv}
\bibliography{gijs}{}

\end{document}